\documentclass[12pt]{amsart}
\usepackage{amssymb}
\usepackage{amsmath}

\allowdisplaybreaks

\newtheorem{thm}{Theorem}[section]
\newtheorem{defn}[thm]{Definition}
\newtheorem{prop}[thm]{Proposition}
\newtheorem{cor}[thm]{Corollary}
\newtheorem{lemma}[thm]{Lemma}
\newtheorem{rema}[thm]{Remark}

\newcommand{\halmos}{\rule{1ex}{1.4ex}}
\newcommand{\pfbox}{\hspace*{\fill}\mbox{$\halmos$}}
\newcommand{\nn}{\nonumber \\}

 \newcommand{\res}{\mbox{\rm Res}}

 \newcommand{\pf}{{\it Proof.}\hspace{2ex}}
 
 \newcommand{\epfv}{\hspace*{\fill}\mbox{$\halmos$}\vspace{1em}}

\newcommand{\wt}{\mbox{\rm wt} \, }

\newcommand{\C}{\mathbb{C}}
\newcommand{\Z}{\mathbb{Z}}

\newcommand{\N}{\mathbb{N}}

        \newcommand{\ba}{\begin{array}}
        \newcommand{\ea}{\end{array}}
        \newcommand{\be}{\begin{equation}}
        \newcommand{\ee}{\end{equation}}
        \newcommand{\bea}{\begin{eqnarray}}
        \newcommand{\eea}{\end{eqnarray}}

\begin{document}
    
\title[Higher level Zhu algebras and modules for VOAs]{Higher level Zhu algebras and modules for vertex operator algebras}
\author{Katrina Barron}
\address{Department of Mathematics, 255 Hurley Hall, University of Notre Dame, Notre Dame, IN 46656}
\email{kbarron@nd.edu}
\author{Nathan Vander Werf}
\address{Department of Mathematics and Statistics, Warner Hall 2036, University of Nebraska, 2504 9th Ave., Kearny, NE 68849} 
\email{vanderserfnp@unk.edu}
\author{Jinwei Yang} 
\address{Department of Mathematics, Yale University, PO Box 208283, New Haven, CT 06520}
\email{jinwei.yang@yale.edu}

\date{October 18, 2018}

\thanks{K. Barron was supported by Simons Foundation Grant 282095}

\begin{abstract}
Motivated by the study of indecomposable, nonsimple modules for a vertex operator algebra $V$, we study the relationship between various types of $V$-modules  and modules for the higher level Zhu algebras for $V$, denoted $A_n(V)$, for $n \in \mathbb{N}$, first introduced by Dong, Li, and Mason in 1998.   We resolve some issues that arise in a few theorems previously presented when these algebras were first introduced, and give examples illustrating the need for certain modifications of the statements of those theorems.   We establish that whether or not $A_{n-1}(V)$ is isomorphic to a direct summand of $A_n(V)$ affects the types of indecomposable $V$-modules which can be constructed by inducing from an $A_n(V)$-module, and in particular whether there are $V$-modules induced from $A_n(V)$-modules that were not already induced by $A_0(V)$.  We give some characterizations of the $V$-modules that can be constructed from such inducings, in particular as regards their singular vectors.  
To illustrate these results, we discuss two examples of $A_1(V)$: when $V$ is the vertex operator algebra associated to either the Heisenberg algebra or the Virasoro algebra.  For these two examples, we show how the structure of $A_1(V)$ in relationship to $A_0(V)$ determines what types of indecomposable $V$-modules can be induced from a module for the level zero versus level one Zhu algebras.  
We construct a family of indecomposable modules for the Virasoro vertex operator algebra that are  logarithmic modules and are not highest weight modules.  
\end{abstract}

\maketitle

\setcounter{equation}{0}
\section{Introduction}

This paper was motivated by a desire to study indecomposable objects in certain module categories for a vertex operator algebra $V$ via modules for the higher level Zhu algebras for $V$, denoted $A_n(V)$, for $n \in \mathbb{N}$.  In \cite{Z}, Zhu introduced an associative algebra, which we denote by $A_0(V)$, for $V$ a vertex operator algebra.  This Zhu algebra has proven to be very useful in understanding the module structure of $V$ for certain classes of modules and certain types of vertex operator algebras.  In particular, in the case that the vertex operator algebra is rational, Frenkel and Zhu, in \cite{FZ}, showed that there is a bijection between the module category for a vertex operator algebra and the module category for the Zhu algebra associated with this vertex operator algebra.  Subsequently, in \cite{DLM}, Dong, Li, and Mason introduced higher level Zhu algebras, $A_n(V)$ for $n \in \mathbb{Z}_+$, and proved many important fundamental results about these algebras.   Dong, Li, and Mason presented several statements that generalize the results of Frenkel and Zhu from the level zero algebras to these higher level algebras, results that mainly focused on the semi-simple setting, e.g., the case of rational vertex operator algebras. 

For an irrational vertex operator algebra, instead of the irreducible modules, indecomposable modules are the fundamental objects in the module category.  To date, it has proven difficult to find examples of vertex operator algebras that have certain nice finiteness properties but non semi-simple representation theory, i.e., so called $C_2$-cofinite irrational vertex operator algebras, which is an important setting for logarithmic conformal field theory.   And in general, in the indecomposable nonsimple setting, the correspondence between the category of such $V$-modules and the category of $A_n(V)$-modules is not well understood.    In this paper we are able to obtain correspondences for certain module subcategories and begin a more systematic study of the settings for which the higher level Zhu algebras become effective tools for understanding and constructing $V$-modules.  

In theory, the higher level Zhu algebras should prove to be important tools for studying indecomposable $\mathbb{N}$-gradable $V$-modules, in particular those that have an increase in the Jordan block size at degree greater than zero with respect to the $\mathbb{N}$-grading.  Whereas the zero level Zhu algebra has been used, for instance in \cite{am1}, \cite{AM},  to determine that important examples, such as the  $\mathcal{W}_p$ triplet vertex operator algebras, do have indecomposable modules, more information about the indecomposables for such vertex operator algebras is given by higher level Zhu algebras, including information necessary to compute the fusion rules for the module category in these non semi-simple settings, e.g. \cite{TW} (see also \cite{NT}).  For $C_2$-cofinite rational vertex operator algebras, Zhu showed in \cite{Z} that the characters (reflective of the $L(0)$ eigenspaces structure)
 of the irreducible modules are closed under the action of modular transformations.  However, as Miyamoto showed in \cite{Miyamoto2004}, for $C_2$-cofinite irrational vertex operator algebras, pseudo-characters for indecomposable non-simple modules must be introduced and considered along with the characters in order to preserve modular invariance, where these so called pseudo-characters detect the generalized eigenspaces of the module.  It is precisely the higher level Zhu algebras which give information about these pseudo-characters.  However despite these important applications of higher level Zhu algebras, so far in practice, higher level Zhu algebras have not been well understood, and techniques for calculating them have not been developed.  This paper is a necessary step toward these goals.

In this paper, we study two functors defined in \cite{DLM}: the functor $\Omega_n/\Omega_{n-1}$   from $\mathbb{N}$-gradable $V$-modules (also called admissible $V$-modules as in \cite{DLM}) to $A_n(V)$-modules; and the functor $L_n$ from $A_n(V)$-modules to $\mathbb{N}$-gradable $V$-modules.  We investigate when the composition of these two functors is isomorphic to the identity morphism in various module categories.  We show that modifications and clarifications are needed for some of the statements in \cite{DLM} to hold, and we give examples to show the necessity of these modifications and clarifications.  We investigate the relationship between indecomposable modules for $A_n(V)$ and certain indecomposable modules for $V$.  We present some sufficient conditions on both $V$-modules and $A_n(V)$-modules for the functors between the restricted module categories to be mutual inverses.  We investigate the question of what types of $V$-modules are constructed from the induction functor $L_n$ and how the structure of $A_n(V)$, in particular as regards $A_{n-1}(V)$, affects the structure of the types of indecomposable modules that can be induced by $L_n$, including for instance when a module can be induced from an $A_n(V)$-module that was not already induced from an $A_{n-1}(V)$-module, and what size Jordan blocks for the $L(0)$ grading operator arise at different levels.  We also give bounds on what degrees singular vectors can reside in for a $V$-module induced by $L_n$ from an $A_n(V)$-module.  

We present $A_1(V)$ for two vertex operator algebras:  the generalized Verma module vertex operator algebras for the Heisenberg and Virasoro algebras, constructed in \cite{BVY-Heisenberg} and \cite{BVY-Virasoro}, respectively.  We describe the nature of the modules for these vertex operator algebras as regards the ring structure of $A_1(V)$, in particular in relationship to $A_0(V)$.    We construct a family of indecomposable nonsimple modules for the Virasoro vertex operator algebra that are logarithmic modules and are not highest weight modules.  We give concrete examples arising from the Virasoro vertex operator algebra motivating the need for the extra conditions we introduce as necessary conditions for the statements of some of the main theorems we prove in this paper to hold.

This paper is organized as follows, in Section 2,  we give basic definitions including the definition of Zhu algebras $A_n(V)$, for $n \in \mathbb{N}$, the functor $\Omega_n$ from the category of $\mathbb{N}$-gradable $V$-modules to the category of $A_n(V)$-modules, and the functor $L_n$ from the category of $A_n(V)$-modules to the category of $\mathbb{N}$-gradable $V$-modules as defined in \cite{DLM}. 

In Section 3,  in the case when $U$ is an $A_n(V)$-module such that {\it no nonzero submodule of} $U$  factors through $A_{n-1}(V)$, we prove in Theorem \ref{mainthm} that $\Omega_n/\Omega_{n-1}(L_n(U)) \cong U$ as $A_n(V)$-modules.  The italics in the previous statement give the extra condition necessary for this statement to hold in contrast to the statement in \cite{DLM}.  We prove Corollary \ref{mainthm-cor} to Theorem \ref{mainthm}  that shows this extra condition is not necessary in the case when $U$ is indecomposable  and $A_n(V)$ decomposes into a direct sum of $A_{n-1}(V)$ and a direct sum complement.  In Section 4, we give examples to illustrate both the need for the extra condition in Theorem \ref{mainthm}  for the case of $n=1$ and the Virasoro vertex operator algebra, and how this extra condition is not needed in, for instance, the case of  $n=1$ for the Heisenberg vertex operator algebra, illustrating how this is predicated on the relationship between $A_1(V)$ and $A_0(V)$ in these two examples.

In Section 3, we show that the functors $L_n$ and  $\Omega_n/\Omega_{n-1}$ are inverse functors when restricted to the categories of simple $A_n(V)$-modules that do not factor through $A_{n-1}(V)$ and simple $\mathbb{N}$-gradable $V$-modules {\it that are generated by their degree $n$ subspace} (or equivalently, {\it that have a nonzero degree $n$ subspace}).  This is a small clarification of Theorem 4.9 of \cite{DLM} where we add the extra condition in italics in the previous statement, for the result to hold.  In Section \ref{Virasoro-n-generated-example}, we give an example involving a simple module for the Virasoro vertex operator algebra to show that this extra condition is necessary for the statement to hold.  

In Section 3, we show that $L_n$ sends indecomposable $A_n(V)$-modules to indecomposable $\mathbb{N}$-gradable $V$-modules.  We go on to identify a certain subcategory of $\mathbb{N}$-gradable $V$-modules on which $L_n \circ \Omega_n/\Omega_{n-1}$ is the identity functor, and we give further restricted categories of $\mathbb{N}$-gradable $V$-modules and $A_n(V)$-modules, respectively, on which the functors restricted to these subcategories are mutual inverses.   More importantly, we identify some characteristics of the types of $V$-modules that can be induced through $L_n$ from an $A_n(V)$-module, for instance properties of the lowest $n$ and first $2n+1$ subspaces and where singular vectors must reside for an $\mathbb{N}$-gradable $V$-module to be induced via $L_n$ from its degree $n$ subspace.

In Section \ref{examples-section}, we present the examples that illustrate the relationship between the level one Zhu algebra and the level zero Zhu algebra and the impact that their ring structures have on their module categories, as well as the necessity of certain additional assumptions for results in Section 3 to hold.  The two examples we consider are the cases for the generalized Verma module vertex operator algebras associated to the Heisenberg and the Virasoro algebras, respectively.  We construct indecomposable $\N$-gradable $V$-modules for each of these vertex operator algebras, $V$, by
presenting the level one Zhu algebra associated with $V$, denoted $A_1(V)$, as constructed in \cite{BVY-Heisenberg} and \cite{BVY-Virasoro}, respectively.  For the Heisenberg vertex operator algebra $A_1(V)$ is isomorphic to $\C[x] \oplus \C[x]$.  Thus the irreducible (indecomposable) modules for $A_1(V)$ are irreducible (indecomposable) modules for either the zero level Zhu algebra $A_0(V)$, which is $\C[\alpha(-1) \mathbf{1}] \cong \C[x]$, or its direct sum component (see Section \ref{Heisenberg-section} for details).  In particular, in the case that $V$ is the Heisenberg vertex operator algebra,
any indecomposable  $A_1(V)$-module $U$ either itself factors through $A_0(V)$ or only its trivial submodule factors through $A_0(V)$.  As a result, there are no new modules induced by $L_1$ from modules for $A_1(V)$ that were not already induced by $L_0$ from $A_0(V)$-modules.  

However, the level one Zhu algebra for the Virasoro vertex operator algebra is isomorphic to $\C[x, y]/( xy)$ as an associative algebra.  Since there is a node at $(0, 0)$ on the curve $xy = 0$,  the level one Zhu algebra provides indecomposable modules that do not factor through the level zero Zhu algebra which is isomorphic to $\C[\omega] \cong \C[x]$, but a nontrivial submodule does factor through $A_0(V)$. For this reason, in this case we can construct explicit indecomposable modules for the Virasoro vertex operator algebra from modules for $A_1(V)$, that are not highest weight modules, i.e. are not induced from an $A_0(V)$-module.  We show that in this case the structure of $A_1(V)$ versus $A_0(V)$ gives rise to a family of modules for each $k \in \mathbb{Z}_+$ that have the property that they are decomposed by the $L(0)$ operator into generalized eigenspaces comprised of Jordan blocks of size $k$ at degree zero and size $k+1$ at degree one and which can not be induced from the level zero Zhu algebra.  

The indecomposable modules for the associative algebra $\C[x, y]/( xy)$ were first studied and classified by Gelfand and Ponomarev \cite{GP} in order to study the representations of Lorentz groups.  A study of indecomposable modules for $\C[x, y]/(xy )$ can be found in \cite{LS}, (cf. \cite{L}, \cite{NR}, \cite{AL}).   In this paper, we present a particular class of modules for $A_1(V)$ for $V$ the Virasoro vertex operator algebra, and in \cite{BVY-Virasoro} we further study the modules of $A_1(V)$ and the resulting modules for the Virasoro vertex operator algebra.  

{\bf Acknowledgements:}  The first author is the recipient of Simons Foundation Collaboration Grant 282095 and greatly appreciates this support.  The authors would like to thank Kiyokazu Nagatomo and the referee for helpful comments on the initial draft of this manuscript.

\section{The notion of level $n$ Zhu algebra $A_n(V)$, and the functors $\Omega_n$, $\Pi_n$, and $L_n$}

Let $V = (V, Y, {\bf 1}, \omega)$ be a vertex operator algebra.   In this section, following \cite{DLM}, we recall the notion of the level $n$ Zhu algebra for $V$, denoted $A_n(V)$, for $n \in \mathbb{N}$, the functor $\Omega_n$ from the category of $\mathbb{N}$-gradable $V$-modules to the category of $A_n(V)$-modules, and the functor $L_n$ from the category of $A_n(V)$-modules to the category of $\mathbb{N}$-gradable $V$-modules. We recall several results of \cite{DLM}, introduce some notation, and prove a lemma that will be useful in proving the main results of this paper.

First, we recall the definition of the level $n$ Zhu algebra $A_n(V)$, for $n \in \N$, introduced in \cite{Z} for $n = 0$, and then generalized to $n >0$ in \cite{DLM}.

\begin{defn}[\cite{Z}, \cite{DLM}]{\rm
For $n \in \N$, let $O_n(V)$ be the subspace of $V$ spanned by elements of the form
\[ u \circ_n v =
\res_x \frac{(1 + x)^{\mathrm{wt}\,  u + n}Y(u, x)v}{x^{2n+2}}
\]
for homogeneous $u \in V$ and for $v \in V$, and by elements of the form $(L(-1) + L(0))v$ for $v \in V$. The vector space $A_n(V)$ is defined to be the quotient space $V/O_n(V)$.

We define the following multiplication on $V$
\[
u *_n v = \sum_{m=0}^n(-1)^m\binom{m+n}{n}\res_x \frac{(1 + x)^{\mathrm{wt} \, u + n}Y(u, x)v}{x^{n+m+1}},
\]
for $v \in V$ and homogeneous $u \in V$, and for general $u \in V$,  $ *_n $  is defined by linearity.  It is shown in \cite{DLM} that with this multiplication, the subspace $O_n(V)$ of $V$ is a two-sided ideal of $V$, and $A_n(V)$ is an associative algebra, called the {\it level $n$ Zhu algebra}.
}
\end{defn}

\begin{rema}\label{epi-remark}
{\em As noted in \cite{DLM}, we have $O_n(V) \subset O_{n-1}(V)$ for $n \in \Z_+$, and thus there is a natural surjective algebra homomorphism from $A_n(V)$ onto $A_{n-1}(V)$ given by $v + O_n(V) \mapsto v+ O_{n-1}(V)$.  If $U$ is a module for $A_n(V)$, then $U$ is said to {\em factor through $A_{n-1}(V)$} if the kernel of the epimorphism from $A_n(V)$ onto $A_{n-1}(V)$, i.e., $O_{n-1}(V)$, acts trivially on $U$, giving $U$ a well-defined $A_{n-1}(V)$-module structure.  }
\end{rema}

Next we recall the definitions of various $V$-module structures.  We assume the reader is familiar with the notion of weak $V$-module for a vertex operator algebra $V$ (cf. \cite{LL}).  

\begin{defn}\label{N-gradable-definition}
{\em An {\it $\mathbb{N}$-gradable weak $V$-module} (also often called an {\it admissible $V$-module} as in \cite{DLM}) $W$ for a vertex operator algebra $V$ is a weak $V$-module that is $\N$-gradable, $W = \coprod_{k \in \N} W(k)$, with $v_m W(k) \subset W(k + \mathrm{wt} v - m -1)$ for homogeneous $v \in V$, $m \in \Z$ and $k \in \N$, and without loss of generality, we can and do assume $W(0) \neq 0$, unless otherwise specified.  We say elements of $W(k)$ have {\it degree} $k \in \mathbb{N}$.

An {\it $\mathbb{N}$-gradable generalized weak $V$-module} $W$ is an $\mathbb{N}$-gradable weak $V$-module that admits a decomposition into generalized eigenspaces via the spectrum of $L(0) = \omega_1$ as follows: $W=\coprod_{\lambda \in{\C}}W_\lambda$ where $W_{\lambda}=\{w\in W \, | \, (L(0) - \lambda \, id_W)^j w= 0 \ \mbox{for some $j \in \mathbb{Z}_+$}\}$, and in addition, $W_{n +\lambda}=0$ for fixed $\lambda$ and for all sufficiently small integers $n$. We say elements of $W_\lambda$ have {\it weight} $\lambda \in \mathbb{C}$.

A {\it generalized $V$-module} $W$ is an $\mathbb{N}$-gradable generalized weak $V$-module where $\dim W_{\lambda}$ is finite for each $\lambda \in \mathbb{C}$.   

An {\it (ordinary) $V$-module} is a generalized $V$-module such that  the generalized eigenspaces $W_{\lambda}$ are in fact eigenspaces, i.e. $W_{\lambda}=\{w\in W \, | \, L(0) w=\lambda w\}$.}
\end{defn}

We will often omit the term ``weak" when referring to $\mathbb{N}$-gradable weak and $\mathbb{N}$-gradable generalized weak $V$-modules.   

The term {\it logarithmic} is also often used in the literature to refer to $\mathbb{N}$-gradable weak generalized modules  or generalized modules. 

\begin{rema}{\em An $\mathbb{N}$-gradable $V$-module with $W(k)$ of finite dimension for each $k \in \mathbb{N}$ is not necessarily a generalized $V$-module since the generalized eigenspaces might not be finite dimensional. }
\end{rema}

We define the {\it generalized graded dimension} of a generalized $V$-module $W = \coprod_{\lambda \in{\C}}W_\lambda$ to be
\begin{equation}
\mathrm{gdim}_q  \, W = q^{-c/24} \,  \sum_{\lambda \in \C}  \, ( \mathrm{dim} \ W_\lambda) \, q^\lambda.
\end{equation}

We recall the functors $\Omega_n$ and $L_n$ for $n \in \N$ defined and studied in \cite{DLM}.  Let $W$ be an $\mathbb{N}$-gradable $V$-module, and let
\begin{equation}
\Omega_n(W) = \{w \in W \; | \; v_iw = 0\;\mbox{if}\; \wt v_i < -n \; 
\mbox{for $v\in V$ of homogeneous weight}\}.
\end{equation}
Then $\Omega_n(W)$ is an $A_n(V)$-module, via the action $[a] \mapsto o(a) = a_{\mathrm{wt} \, a -1}$ for $a \in V$ and $[a] = a + O_n(V)$.

\begin{rema} {\em The functor $\Omega_n$ from $\N$-gradable $V$-modules to $A_n(V)$-modules we have defined here is the same functor defined in \cite{DLM}, but is not the functor called $\Omega_n$ in \cite{V}; rather in \cite{V} the functor denoted $\Omega_n$ is just the projection functor onto the $n$th graded subspace of $W$.  We will denote this projection functor by $\Pi_n$; that is for $W$ an $\mathbb{N}$-gradable $V$-module, we have $\Pi_n(W) = W(n)$.  }
\end{rema}

In order to define the functor $L_n$ from the category of $A_n(V)$-modules to the category of $\mathbb{N}$-gradable $V$-modules, we need several notions, including the notion of the universal enveloping algebra of $V$, which we now define.  

Let
\begin{equation}
\hat{V} = \C[t, t^{-1}]\otimes V/D\C[t, t^{-1}]\otimes V,
\end{equation}
where $D = \frac{d}{dt}\otimes 1 + 1 \otimes L(-1)$. For $v \in V$, let $v(m) = v \otimes t^m +  D\C[t, t^{-1}]\otimes V \in \hat{V}$.  Then $\hat{V}$ can be given the structure of a $\Z$-graded Lie algebra as follows:  Define the degree of $v(m)$ to be $\wt v - m - 1$ for homogeneous $v \in V$, and define the Lie bracket on $\hat{V}$ by
\begin{equation}\label{bracket}
[u(j), v(k)] = \sum_{i = 0}^{\infty}\binom{j}{i}(u_iv)(j+k-i),
\end{equation}
for $u, v \in V$, $j,k \in \Z$.
Denote the homogeneous subspace of degree $m$ by $\hat{V}(m)$. In particular, the degree $0$ space of $\hat{V}$, denoted by $\hat{V}(0)$, is a Lie subalgebra.

Denote by $\mathcal{U}(\hat{V})$ the universal enveloping algebra of the Lie algebra $\hat{V}$.  Then $\mathcal{U}(\hat{V})$ has a natural $\Z$-grading induced from $\hat{V}$, and we denote by $\mathcal{U}(\hat{V})_k$ the degree $k$ space with respect to this grading, for $k \in \Z$.

Given a weak $V$-module $W$, consider the following linear map
\begin{eqnarray}\label{defining-phi}
\varphi_W : \ \   \ \ \mathcal{U}(\hat{V}) &\longrightarrow & \mathrm{End} (W)\\
v_1(m_1)v_2(m_2) \cdots v_k(m_k) & \mapsto & (v_1^W)_{m_1}(v_2^W)_{m_2} \cdots (v_k^W)_{m_k}  ,\nonumber
\end{eqnarray}
where $v^W_m$ is the coefficient of $x^{-m-1}$ in the vertex operator action of $Y_W(v,x)$ on the module $W$.  We will often denote $\varphi_W$ by $\varphi$ and $v^W_m$ by $v_m$  if the module is $V$ itself or if $W$ is clearly implied.

We shall need the following lemma in the proof of our main theorem:
\begin{lemma}\label{l1}
We have
\begin{equation}\label{split}
o(O_n(V)) \subseteq \coprod_{i > n} \varphi(\mathcal{U}(\hat{V})_i) \varphi (\mathcal{U}(\hat{V})_{-i}).
\end{equation}
\end{lemma}

\pf  By the $L(-1)$-derivative property in $V$, we have that $o(L(-1)v) = -o(L(0)v)$, and thus $o((L(-1) + L(0))v) = 0$, showing (\ref{split}) holds trivially for the elements of the form $(L(-1)+L(0))v$ in $O_n(V)$.

Next, for homogeneous $u, v \in V$, we have
\begin{eqnarray}
o \left(\res_x \frac{(1 + x)^{\mathrm{wt} \, u+n}Y(u, x)v}{x^{2n+2}} \right) \nonumber
&=& \sum_{j \in \mathbb{N}}  \binom{\mathrm{wt} \, u + n}{j} o(u_{j-2n-2}v)\nn
&= & \sum_{j \in \mathbb{N}} \binom{\mathrm{wt} \,  u + n}{j}\res_{x_2}\res_{x_1-x_2}  Y(Y(u, x_1-x_2)v, x_2)\nn
&& \quad (x_1-x_2)^{j-2n-2}x_2^{\mathrm{wt} \, u + \mathrm{wt} \,  v -j+2n}\nn
&= & \res_{x_2}\res_{x_1-x_2} Y(Y(u, x_1-x_2)v, x_2)\frac{x_1^{\mathrm{wt} \,  u + n} x_2^{\mathrm{wt} \, v+n}}{(x_1-x_2)^{2n+2}}\nn
&= & \res_{x_1}\res_{x_2} Y(u, x_1)Y(v, x_2)\frac{x_1^{\mathrm{wt} \, u + n}x_2^{\mathrm{wt} \,  v+n}}{(x_1-x_2)^{2n+2}}\nn
&& - \res_{x_1}\res_{x_2} Y(v, x_2)Y(u, x_1)\frac{x_1^{\mathrm{wt} \,  u + n}x_2^{\mathrm{wt}\, v+n}}{(-x_2+x_1)^{2n+2}}, \label{lemma-proof}
\end{eqnarray}
where the last equality follows from the Jacobi identity on $V$.

The first term on the right hand side of (\ref{lemma-proof}) satisfies
\begin{eqnarray*}
\lefteqn{\res_{x_1}\res_{x_2} Y(u, x_1)Y(v, x_2)\frac{x_1^{\mathrm{wt} \, u + n}x_2^{\mathrm{wt} \, v+n}}{(x_1-x_2)^{2n+2}}}\nonumber \\
&= &\sum_{j \in \mathbb{N}}(-1)^{j}\binom{-2n-2}{j}\res_{x_1}\res_{x_2} Y(u, x_1)Y(v, x_2) x_1^{\mathrm{wt} \, u -n-2-j}x_2^{\mathrm{wt} \, v+n+j}\nn
&= & \sum_{j \in \mathbb{N}}(-1)^{j}\binom{-2n-2}{j} u_{\mathrm{wt} \, u -n-2-j}v_{\mathrm{wt} \, v+n+j}\nn
&=& \sum_{k\leq -n-1} (-1)^{n+k+1} \binom{-2n-2}{-n-k-1} u_{\mathrm{wt} \,  u +k -1}v_{\mathrm{wt} \, v -k -1}\nn
&\in & \coprod_{i > n}\varphi(\mathcal{U}(\hat{V})_i) \varphi(\mathcal{U}(\hat{V})_{-i}).
\end{eqnarray*}

Similarly, the second term of the right hand side of (\ref{lemma-proof}) satisfies 
\begin{eqnarray*}
\lefteqn{\res_{x_1}\res_{x_2} Y(v, x_2)Y(u, x_1)\frac{x_1^{\mathrm{wt} \, u + n}x_2^{\mathrm{wt} \, v+n}}{(-x_2 + x_1)^{2n+2}}}\nonumber \\
&= &\sum_{j \in \mathbb{N}}(-1)^{j}\binom{-2n-2}{j}\res_{x_1}\res_{x_2} Y(v, x_2)Y(u, x_1) x_1^{\mathrm{wt} \, u +n + j}x_2^{\mathrm{wt} \, v-n-j-2}\nn
&= & \sum_{j \in \mathbb{N}}(-1)^{j}\binom{-2n-2}{j} v_{\mathrm{wt} \, v -n-2-j}u_{\mathrm{wt} \, u+n+j}\nn
&=& \sum_{k\leq -n-1} (-1)^{n+k+1} \binom{-2n-2}{-n-k-1} u_{\mathrm{wt}\,  v +k -1}v_{\mathrm{wt}\,  u -k -1}\nn
&\in & \coprod_{i > n}\varphi(\mathcal{U}(\hat{V})_i) \varphi(\mathcal{U}(\hat{V})_{-i}).
\end{eqnarray*}
\pfbox\\

We can regard $A_n(V)$ as a Lie algebra via the bracket $[u,v] = u *_n v - v *_n u$, and then the map $v( \mathrm{wt} \, v -1) \mapsto v + O_n(V)$ is a well-defined Lie algebra epimorphism from $\hat{V}(0)$ onto $A_n(V)$.

Let $U$ be an $A_n(V)$-module.  Since $A_n(V)$ is naturally a Lie algebra homomorphic image of $\hat{V}(0)$, we can lift $U$ to a module for the Lie algebra $\hat{V}(0)$, and then to a module for $P_n = \bigoplus_{p > n}\hat{V}(-p) \oplus \hat{V}(0) = \bigoplus_{p < -n} \hat{V}(p) \oplus \hat{V}(0)$ by letting $\hat{V}(-p)$ act trivially for $p\neq 0$.  Define
\[
M_n(U) = \mbox{Ind}_{P_n}^{\hat{V}}(U) = \mathcal{U}(\hat{V})\otimes_{\mathcal{U}(P_n)}U.
\]

We impose a grading on $M_n(U)$ with respect to $U$, $n$, and the $\Z$-grading on $\mathcal{U}(\hat{V})$, by letting $U$ be degree $n$, and letting $M_n(U)(k)$, for $k \in  \Z$, to be the subspace of $M_n(U)$ induced from $\hat{V}$, i.e., $M_n(U)(k) = \mathcal{U}(\hat{V})_{k-n}U$. 

For $v \in V$, define $Y_{M_n(U)}(v,x) \in (\mathrm{End} (M_n(U)))((x))$ by
\begin{equation}\label{define-Y_M}
Y_{M_n(U)}(v,x) = \sum_{m\in\mathbb{Z}} v(m) x^{-m-1}.
\end{equation} 

Let $W_{A}$ be the subspace of $M_n(U)$ spanned linearly by the coefficients of 
\begin{multline}\label{relations-for-M}
(x_0 + x_2)^{\mathrm{wt} \, v + n} Y_{M_n(U)}(v, x_0 + x_2) Y_{M_n(U)}(w, x_2) u \\ 
- (x_2 + x_0)^{\mathrm{wt} \, v + n} Y_{M_n(U)}(Y(v, x_0)w, x_2) u
\end{multline}
for $v,w \in V$, with $v$ homogeneous, and $u \in U$.  Set
\[ \overline{M}_n(U) = M_n(U)/\mathcal{U} (\hat{V})W_A .\]

It is shown in \cite{DLM} that if $U$ is an $A_n(V)$-module that does not factor through $A_{n-1}(V)$, then $\overline{M}_n(U) = \bigoplus_{k \in \mathbb{N}} \overline{M}_n(U) (k)$ is an $\mathbb{N}$-gradable $V$-module with $\overline{M}_n(U) (0)\neq 0$ and $\overline{M}_n(U) (n) \cong U$ as an $A_n(V)$-module. Note that the condition that $U$ itself does not factor though $A_{n-1}(V)$ is indeed a necessary and sufficient condition for $\overline{M}_n(U) (0)\neq 0$ to hold.  

It is also observed in \cite{DLM} that $\overline{M}_n(U)$ satisfies the following universal property:  For any weak $V$-module $M$ and any $A_n(V)$-module homomorphism $\phi: U \longrightarrow \Omega_n(M)$, there exists a unique weak $V$-module homomorphism $\Phi: \overline{M}_n(U) \longrightarrow M$, such that $\Phi \circ \iota = \phi$ where $\iota$ is the natural injection of $U$ into $\overline{M}_n(U)$. This follows from the fact that $\overline{M}_n(U)$ is generated by $U$ as a weak $V$-module, again with the possible need of a grading shift.

Let $U^* = \mbox{Hom}(U, \C)$.  As in the construction in \cite{DLM}, we can extend $U^*$ to $M_n(U)$ by first an induction to $M_n(U)(n)$ and then by letting $U^*$ annihilate $\bigoplus_{k \neq n} M_n(U)(k)$.  In particular, we have that elements of $M_n(U)(n) = \mathcal{U}(\hat{V})_0U$ are spanned by elements of the form 
\[o_{p_1}(a_1) \cdots o_{p_s}(a_s)U\]
where $s \in \mathbb{N}$, $p_1 \geq \cdots \geq p_s$, $p_1 + \cdots + p_s =0$, $p_i \neq 0$, $p_s \geq -n$, $a_i \in V$ and $o_{p_i}(a_i) = (a_i)(\mathrm{wt} \, a_i - 1 - p_i)$. Then inducting on $s$ by using Remark 3.3 in \cite{DLM} to reduce from length $s$ vectors to length $s-1$ vectors, we have a well-defined action of $U^*$ on $M_n(U)(n)$.  

Set
\[
J = \{v \in M_n(U) \, | \, \langle u', xv\rangle = 0 \;\mbox{for all}\; u' \in U^{*}, x \in \mathcal{U}(\hat{V})\}
\]
and
\[
L_n(U) = M_n(U)/J.
\]

\begin{rema}\label{L-a-V-module-remark} {\em It is shown in \cite{DLM}, Propositions 4.3, 4.6 and 4.7,  that if $U$ does not factor through $A_{n-1}(V)$, then $L_n(U)$ is a well-defined $\mathbb{N}$-gradable $V$-module with $L_n(U)(0) \neq 0$; in particular, it is shown that $\mathcal{U}(\hat{V})W_A \subset J$, for $W_A$ the subspace of $M_n(U)$ spanned by the coefficients of (\ref{relations-for-M}), i.e., giving the associativity relations for the weak vertex operators on $M_n(U)$.}
\end{rema}

\section{Main results}

We have the following theorem which is a necessary modification to what was presented as Theorem 4.2 in \cite{DLM}:

\begin{thm}\label{mainthm}
For $n \in \N$, let $U$ be a nonzero $A_n(V)$-module such that if $n>0$, then $U$ does not factor through $A_{n-1}(V)$. Then $L_n(U)$ is an $\mathbb{N}$-gradable $V$-module with $L_n(U)(0) \neq 0$.  If we assume further that there is no nonzero submodule of $U$ that factors through $A_{n-1}(V)$, then $\Omega_n/\Omega_{n-1}(L_n(U)) \cong U$.
\end{thm}

\pf  For $n \in \N$, let $U$ be a nonzero $A_n(V)$-module such that if $n>0$, then $U$ does not factor through $A_{n-1}(V)$.  By Remark \ref{L-a-V-module-remark}, $L_n(U)$ is a well-defined $\mathbb{N}$-gradable $V$-module with $L_n(U)(0) \neq 0$, and 
\begin{equation}\label{in-J}
\mathcal{U}(\hat{V})W_A \subset J.
\end{equation}

We have $U \subset M_n(U)(n) \subset M_n(U)$, and for $u \in U$, we have $\langle U^*, u \rangle = 0$ if and only if $u = 0$.   Therefore $J \cap U = 0$, implying that $U$ is naturally isomorphic to a subspace of $L_n(U)$, namely $U + J$.  We denote $u + J$ by $\bar{u}$ for $u \in U$.

Since $U$ is an $A_n(V)$-module, we have that $o(v) \cdot U = 0$ for all $v \in O_n(V)$, and thus $\bar{u} \in \Omega_n(L_n(U))$ for all $u \in U$.  But 
\[(U + J)\cap \Omega_{n-1}(L_n(U)) = \{u + J \; | \; u \in U \ \mathrm{and} \ v_i (u +J) = 0 \ \mathrm{if} \ \mathrm{wt} \, v_i < -n + 1 \}\] 
is annihilated by $\mathcal{U}(\hat{V})_{-i}$ for $i < -n+1$ and therefore, by Lemma \ref{l1},  is annihilated by the action of $O_{n-1}(V)$, implying that $(U + J) \cap \Omega_{n-1}(L_n(U))$ is a nontrivial $A_n(V)$-submodule of $U +J$ that factors through $A_{n-1}(V)$.  

Now assume that $U$ has no nonzero submodule that factors through $A_{n-1}(V)$.  This and the fact that $U\cap J = 0$, implies that $(U + J)\cap \Omega_{n-1}(L_n(U)) = 0$.  Therefore there is an injection of $A_n(V)$-modules $\iota: U \hookrightarrow \Omega_n/\Omega_{n-1} (L_n(U))$ given by $u \mapsto \bar{u} + \Omega_{n-1}(L_n(U))$.  In particular, this implies that $L_n(U) (0) \neq 0$ since $U$ is assumed to be nonzero.

Next we show the injection $\iota$ is surjective, i.e., if $\bar{w} + \Omega_{n-1}(L_n(U)) \in  \Omega_n/\Omega_{n-1} (L_n(U))$, for $\bar{w} \in \Omega_n(L_n(U))$, then there exists $u \in U$ such that $\bar{u} = \bar{w} + v$ for $v \in \Omega_{n-1} (L_n(U))$, i.e., $\iota(u) = \bar{w} + \Omega_{n-1} (L_n(U))$.  For convenience, we denote the coset $\bar{w} + \Omega_{n-1} (L_n(U))$ by $\hat{w}$.

Let $\bar{w} \in \Omega_n(L_n(U))$, and let $w \in M_n(U)$ be its preimage under the canonical projection.   

If $\mathrm{deg} \, w < n$, then for $v \in V$ with $\wt v_i < -n + 1$, we have $\wt v_i w < 0$ and thus $\bar{w} \in \Omega_{n-1}(L_n(U))$, and $\iota(0) = \hat{0} = \hat{w}$.   

For the case when $\mathrm{deg} \, w = n$, by (\ref{in-J}) there is a natural surjection 
\begin{eqnarray}
\pi : \overline{M}_n(U) = M_n(U)/ \mathcal{U}(\hat{V})W_A &\longrightarrow & L_n(U) = M_n(U)/J \\
w + \mathcal{U}(\hat{V})W_A & \mapsto & w +J. \nonumber 
\end{eqnarray}
This coupled with the fact that $\overline{M}_n(U)(n) \cong U$ implies that if $\mathrm{deg} \, w = n$, then there exists a $u \in U$ such that $u +  \mathcal{U}(\hat{V})W_A = w +  \mathcal{U}(\hat{V})W_A \stackrel{\pi}{\mapsto} w + J  = \bar{w} \in \Omega_n(L_n(U))$, proving that $\iota(u) = \bar{w} + \Omega_{n-1}(L_n(U))$.

Finally if $\mathrm{deg} \, w >n$, we will show that $\hat{w} = \hat{0}$ and thus $\iota(0) = \hat{w}$.  Suppose not.  If $\hat{w} \neq \hat{0}$, then in particular $w \notin J$, and thus there exists $x \in \mathcal{U}(\hat{V})_{n - \mathrm{deg}\,  w}$ such that $\langle u', xw \rangle \neq 0$ for some $u' \in U^*$.  Consider the $A_n(V)$-module generated by $xw \in M_n(U)$, defined as follows:
For $v + O_n(V) \in A_n(V)$ and $u \in U$, let $(v + O_n(V)) \cdot u = v(\mathrm{wt} \, v - 1) u$.  Under this action, let $\overline{U} = A_n(V) \cdot xw$.  By Lemma \ref{l1}, $O_n(V) \cdot \overline{U} = \coprod_{i >n} \varphi(\mathcal{U}(\hat{V})_i )\varphi(\mathcal{U}(\hat{V})_{-i}) \cdot\overline{U} =  \coprod_{i >n} \mathcal{U}(\hat{V})_i )\mathcal{U}(\hat{V})_{-i} \overline{U}  = 0$.  However this implies that $O_{n-1}(V) \cdot \overline{U} = \coprod_{i >n-1} \mathcal{U}(\hat{V})_i \mathcal{U}(\hat{V})_{-i} \overline{U} = \mathcal{U}(\hat{V})_n \mathcal{U}(\hat{V})_{-n}  \overline{U}$. But $\hat{w} \neq \hat{0}$ in $\Omega_n/\Omega_{n-1} (L_n(U))$ also implies, in particular,  that $v_i (w + J) = 0$ if $\wt v_i < -n$.  Noting then that 
\[\mathcal{U}(\hat{V})_{-n} \cdot \overline{U} = \{ v_i o(v) x w \; | \;  v_i \in \mathcal{U}(\hat{V})_{-n}, \,  v \in A_n(V) \, \mathrm{with} \, \mathrm{wt} \, o(v) = 0, \, \mathrm{and \ wt} \, x = n- \mathrm{wt} \, w \},\]
and thus we have that $\wt v_i o(v) x = - \wt w < -n$, we then must have that $\mathcal{U}(\hat{V})_{-n} \cdot \overline{U} = 0$. Therefore
$o(O_{n-1}(V)) \cdot \overline{U} = \mathcal{U}(\hat{V})_n \mathcal{U}(\hat{V})_{-n} \cdot \overline{U} = 0$, implying that $\bar{U} = A_n(V) \cdot xw$ is a submodule of $U$ that  factors through $A_{n-1}(V)$.  By assumption then, we must have that $xw = 0$, and thus $w \in J$, implying that $\iota(0) = \hat{0} = \hat{w}$.  
\epfv

Note that it is trivially true that $\Pi_n (L_n(U)) \cong U$ as was observed in \cite{V} (where $\Pi_n$ in \cite{V} is denoted by $\Omega_n$).

Theorem 4.2 in \cite{DLM} only imposes the condition on $U$ that it be an $A_n(V)$-module that itself does not factor through $A_{n-1}(V)$ in order for $\Omega_n/\Omega_{n-1}(L_n(U)) \cong U$ to hold,
rather than the condition that no nonzero submodule of $U$ factor through $A_{n-1}(V)$.  In Section \ref{Virasoro-factor-through-example}, we give an example that shows that this extra condition is necessary.   This example is based on a module for the Virasoro vertex operator algebra and the level one Zhu algebra.  In Section \ref{Heisenberg-section}, we show that this extra condition is not needed at level $n=1$ for the Heisenberg vertex operator algebra.  Furthermore in these examples, we observe why these cases are different in regards to the structure of $A_1(V)$ for the Heisenberg versus the Virasoro vertex operator algebras as regards the algebra $A_0(V)$.

One of the main reasons we are interested in Theorem \ref{mainthm} is what it implies for the question of when modules for the higher level Zhu algebras give rise to indecomposable nonsimple modules for $V$.   In \cite{V} it is claimed that if $A_n(V)$ is a finite dimensional semisimple algebra for all $n \in \mathbb{N}$, then $V$ is rational.  But we are interested in the irrational case and in particular when $A_n(V)$ can be used to construct certain indecomposable nonsimple modules for $V$ of interest in the irrational setting, particularly in the $C_2$-cofinite case.

To this end, we have the following two corollaries to Theorem \ref{mainthm}:

\begin{cor}\label{mainthm-first-cor}  Suppose that for some fixed $n \in \mathbb{Z}_+$, $A_n(V)$ has a direct sum decomposition $A_n(V) \cong A_{n-1}(V) \oplus A'_n(V)$, for $A'_n(V)$ a direct sum complement to $A_{n-1}(V)$, and let $U$ be an $A_n(V)$-module.   If $U$ is trivial as an  $A_{n-1}(V)$-module, then   $\Omega_n/\Omega_{n-1} (L_n(U) ) \cong U$. 
\end{cor}

\pf  If $A_n(V) \cong A_{n-1}(V) \oplus A'_n(V)$, then any $A_n(V)$-module $U$ decomposes into $U = U_{n-1} \oplus U'$ where $U_{n-1}$ is an $A_{n-1}(V)$-module and $U'$ is an $A'_n(V)$-module.  If $U$ has no nontrivial submodule that factors through $A_{n-1}(V)$ which is true if and only if $U_{n-1} = 0$, i.e., if and only if $U$ is trivial as an $A_{n-1}(V)$-module, then Theorem \ref{mainthm} implies $\Omega_n/\Omega_{n-1} (L_n(U) ) \cong U$.
\epfv

An example of this setting is given in Section 4, below, namely that of the Heisenberg vertex operator algebra and the level one Zhu algebra. 

\begin{cor}\label{mainthm-cor}
Let $n \in \mathbb{Z}_+$ be fixed, and let $U$ be a nonzero indecomposable $A_n(V)$-module such that there is no nonzero submodule of $U$ that can factor through $A_{n-1}(V)$. Then $L_n(U)$ is a nonzero indecomposable $\mathbb{N}$-gradable $V$-module generated by its degree $n$ subspace, $L_n(U) (n) \cong U$, and satisfying 
\[\Omega_n/\Omega_{n-1} (L_n(U)) \cong L_n(U)(n) \cong U\] 
as $A_n(V)$-modules.

Furthermore if $U$ is a simple $A_n(V)$-module, then $L_n(U)$ is a simple $V$-module as well.
\end{cor}

\pf Suppose
\[
L_n(U) = W_1 \oplus W_2,
\]
where $W_1$ and $W_2$ are nonzero $\mathbb{N}$-gradable $V$-modules. Then by Theorem \ref{mainthm}, and linearity, we have that
\begin{eqnarray*}
U & \cong & \Omega_n/\Omega_{n-1}(L_n(U)) \\
&=& \Omega_n/\Omega_{n-1}(W_1 \oplus W_2)\\
&=& \Omega_n/\Omega_{n-1}(W_1) \oplus \Omega_n/\Omega_{n-1}(W_2).
\end{eqnarray*}
Since $U$ is indecomposable, $\Omega_n/\Omega_{n-1}(W_i) = 0$ for $i = 1$ or $2$. Without loss of generality assume that $\Omega_n/\Omega_{n-1}(W_1) = 0$.  Then in 
$L_n(U)$, we have that  $0 \neq (U + J) \cap W_1(n) \subset (U +J) \cap \Omega_n(W_1) = (U + J) \cap \Omega_{n-1}(W_1)$, which is an $A_{n-1}(V)$-module.  Thus $\{ u \in U \; | \; u + J \in W_1(n)\} \subset U$ is a nonzero submodule of $U$ that factors through $A_{n-1}(V)$,  contradicting our assumption.  Thus it follows that either $W_1$ or $W_2$ is zero, and $L_n(U)$ is indecomposable.

It is obvious that $L_n(U)$ is generated by $L_n(U) (n) \cong U$ and the fact that this is isomorphic as an $A_n(V)$-module to  $\Omega_n/\Omega_{n-1}(L_n(U))$ follows directly from Theorem \ref{mainthm}. 

Finally we observe that the proof of Lemma 4.8 in \cite{DLM} holds, proving the last statement.  
\epfv

\begin{defn}\label{catasso} 
For $n \in \mathbb{N}$, denote by $\mathcal{A}_{n,n-1}$ the category of $A_n(V)$-modules such that there are no nonzero submodules that can factor through $A_{n-1}(V)$.
\end{defn}

\begin{defn}\label{V-category-definition} 
For $n \in \N$, denote by $\mathcal{V}_n$ the category of weak $V$-modules whose objects $W$ satisfy:  $W$ is $\mathbb{N}$-gradable with $W(0) \neq 0$; and $W$ is generated by $W(n)$.
\end{defn}



With these definitions, we have that Theorem \ref{mainthm} shows that $\Omega_n/\Omega_{n-1} \circ L_n$ is the identity functor on the category $\mathcal{A}_{n,n-1}(V)$, and Corollary \ref{mainthm-cor} shows that $L_n$ sends indecomposable objects in $\mathcal{A}_{n,n-1}$ to indecomposable objects in $\mathcal{V}_n$.  Natural questions that arise then are: What subcategory of $V$-modules does one need to restrict to so that these functors are inverses of each other on this restricted category?  But more importantly, what can be said about the correspondences between the subcategories of simple versus indecomposable objects?   Below we further investigate these questions.

\begin{rema}{\em
Theorem 4.9 of \cite{DLM} states that the functors $L_n$ and $\Omega_n/\Omega_{n-1}$ induce mutually inverse bijections between the isomorphism classes of simple objects in the category $\mathcal{A}_{n,n-1}$ and the isomorphism classes of simple objects in the category of $\mathbb{N}$-gradable $V$-modules. However it is necessary to add the condition that the degree $n$ subspace be nonzero for this to hold, as our example in Section \ref{Virasoro-n-generated-example} shows.  That is, even in the case of simple objects, if the $\mathbb{N}$-gradable $V$-module is not generated by its degree $n$ subspace (e.g., if the degree $n$ subspace is zero), then the bijection can fail.  With this small and obvious addition, we have the theorem below---a clarification of Theorem 4.9 in \cite{DLM}; see also \cite{V}. }
\end{rema}

\begin{thm}\label{simple-theorem}
$L_n$ and $\Omega_n/\Omega_{n-1}$ are equivalences when restricted to the full subcategories of completely reducible $A_n(V)$-modules whose irreducible components do not factor through $A_{n-1}(V)$ and completely reducible $\mathbb{N}$-gradable $V$-modules that are generated by their degree $n$ subspace (or equivalently, that have a nonzero degree $n$ subspace), respectively.  In particular, $L_n$ and $\Omega_n/\Omega_{n-1}$ induce naturally inverse bijections on the isomorphism classes of simple objects in the category $\mathcal{A}_{n, n-1}$ and isomorphism classes of simple objects in $\mathcal{V}_n$.  
\end{thm}

\begin{rema}
{\em We note that Theorem 4.10 of \cite{DLM} in the case of $V$ rational (i.e., in the semisimple setting) also must be modified accordingly with the extra condition that the category of $\mathbb{N}$-gradable $V$ modules must be restricted to the subcategory of modules generated by their degree $n$ subspace in order for the statement of the theorem to hold; see also \cite{V}, where in the rational setting $\Omega_n/\Omega_{n-1}$ can be replace by $\Pi_n$ in many of the statements of \cite{DLM}.  }
\end{rema}

We have the following result which gives a stronger statement than Corollary \ref{mainthm-cor} and a way of constructing indecomposable $\mathbb{N}$-gradable $V$-modules, which we employ in Section \ref{Virasoro-factor-through-example}:

\begin{prop} Let $U$ be an indecomposable $A_n(V)$-module that does not factor through $A_{n-1}(V)$. Then $L_n(U)$ is an indecomposable $\mathbb{N}$-gradable $V$-module generated by its degree $n$ subspace. 

Furthermore, if $U$ is finite dimensional, then $L_n(U)$ is an indecomposable $\mathbb{N}$-gradable generalized $V$-module.
\end{prop}

\pf  Suppose $L_n(U) = W_1 \oplus W_2$, where $W_1$ and $W_2$ are nonzero $\N$-gradable submodules of $L_n(U)$. Then $U +J  = L_n(U)(n) = W_1(n) \oplus W_2(n)$.
Since $U$ is indecomposable, we have that either $U +J = W_1(n)$ or $U +J = W_2(n)$. Without loss of generality, assume that $U +J = W_2(n)$. Then $ W_1 \cap (U +J) = 0$. Let $\widetilde{W_1}$ be the preimage of $W_1$ in $M_n(U)$. Then $\widetilde{W_1} \cap U = 0$, and hence $\widetilde{W_1}(n)\subset J$. Then since $L_n(U)$ is generated by $L_n(U)(n)$, we have that $W_1$ is generated by $W_1(n)$, and thus we have $\widetilde{W}_1 \subset J$. Therefore $W_1 = \widetilde{W_1}/J = 0$, contradicting the assumption that $W_1$ was nonzero, and proving that $L_n(U)$ is indecomposable.   

Now assume that $U$ is finite dimensional.  Since $L(0) = o(\omega)$ preserves $U$, and $U$ is finite dimensional, we have that $U +J = L_n(U)(n)$ can be decomposed into a direct sum of generalized eigenspaces for $L(0)$. However, since $\omega$ is in the center of $A_n(V)$, the distinct generalized eigenspaces of $U$ with respect to $L(0)$ are distinct $A_n(V)$-submodules of $U$.  Therefore there exists $\lambda \in \mathbb{C}$ and $j \in \mathbb{Z}_+$ such that in $L_n(U)$, we have $(L(0) - \lambda \, id_{U+J})^j (U+J)= 0$. 

Then $L_n(U)$ has a $\mathbb{C}$-grading with respect to the eigenvalues of $L(0)$ induced by $M_n(U)$ and the eigenvalue $\lambda$ of $L(0)$ on $U$, given by $M_n(U)(k) = M_n(U)_{\lambda -n + k}$, proving that $L_n(U)$ is an $\mathbb{N}$-gradable generalized $V$-module.
\epfv

In general $\Omega_n/\Omega_{n-1}$ will not send indecomposable objects in $\mathcal{V}_{n}$ to indecomposable $A_n$-modules, or for that matter will $\Pi_n$.  And so we have the following questions:   If $W$ is an indecomposable object in $\mathcal{V}_{n}$, when are $W(n) = \Pi_n W$ and $\Omega_n/\Omega_{n-1}(W)$ indecomposable  $A_n(V)$-modules? Furthermore, and more importantly for our purposes, what types of indecomposable modules can be constructed from the functor $L_n$?  We begin to answer some of these questions below.

It is easy to see the following:
\begin{prop}
Let $W$ be an indecomposable object in $\mathcal{V}_n$. Then the $A_n(V)$-module $W(n)$ cannot be decomposed into a direct sum of subspaces $U_1$ and $U_2$ such that $\langle U_1\rangle \cap \langle U_2\rangle = 0$, where $\langle U_i\rangle$ denotes the $V$-submodule of $W$ generated by $U_i$, for $i = 1, 2$.
\end{prop}
\pf Otherwise, assume $W(n) = U_1 \oplus U_2$ such that $\langle U_1\rangle \cap \langle U_2\rangle = 0$. Since $W$ is generated by $W(n)$, we have $W = \langle U_1\rangle \oplus \langle U_2\rangle$, contradicting the assumption that $W$ was indecomposable. \epfv

Regarding the question:  When is $L_n( \Omega_n/\Omega_{n-1} (W)) \cong W$ for $W$ an $\mathbb{N}$-gradable $V$-module? It is clear that we must at least have that $W$ is an object in the category $\mathcal{V}_n$, i.e., $W$ must be generated by its degree $n$ subspace $W(n)$.  However, in general we have that 
\[ \bigoplus_{k = 0}^n W(k) \subset \Omega_n(W) \]
with equality holding if, for instance, $W$ is simple.  In Section \ref{Virasoro-factor-through-example}, we give an example, that of the Virasoro vertex operator algebra, to show that in the indecomposable case equality will not necessarily hold.  But we do have the following sufficient criteria,  where below we denote the cyclic submodule of $W$ generated by $w\in W$ by $Vw$:

\begin{thm}\label{last-thm}
Let $W$ be an $\mathbb{N}$-gradable $V$-module that is generated by $W(n)$ such that $\Omega_j(W) = \bigoplus_{k=0}^j W(k)$, for $j = n$ and $n-1$.   Then $L_n(\Omega_n/\Omega_{n-1}(W))$ is naturally isomorphic to a quotient of $W$. 

Furthermore, suppose $W$ also satisfies the property that for any $w \in W$, $w =0$ if and only if $V w \cap W(n)  = 0$. Then  
\[L_n(\Omega_n/\Omega_{n-1}(W) ) \cong W.\]
\end{thm}

\pf 
Consider the $A_n(V)$-module injection from $W(n)$ into $\Omega_n(W) = \bigoplus_{k = 0}^n W(k)$.  Then by the universal property of $\overline{M}_n(W(n))$, and since $W$ is generated by $W(n)$, there exists a unique $V$-module surjection 
\begin{eqnarray*}
\Phi : \overline{M}_n(W(n)) &\longrightarrow& W\\
xw + \mathcal{U}(\hat{V})W_A &\mapsto& \varphi(x) w
\end{eqnarray*}
for $x \in \mathcal{U}(\hat{V})$ and $w \in W(n)$, where $\varphi$ is defined in (\ref{defining-phi}). 

Letting $\overline{J} = J/\mathcal{U}(\hat{V})W_A$, where 
\[J = \{ w \in M_n(W(n)) \, | \,  \langle u', xw \rangle = 0 \mbox{ for all $u' \in W(n)^*$, and $x \in \mathcal{U}(\hat{V})$} \}, \]
we have that $\overline{J}$ is naturally an $\mathbb{N}$-gradable $V$-submodule of $\overline{M}_n(W(n))$.  Then $\mathrm{ker} \, \Phi \subset \overline{J}$, and thus 
\[L_n(\Omega_n/\Omega_{n-1} (W)) \cong L_n(W(n)) = M_n(W(n))/J \cong \overline{M}_n(W(n))/\overline{J}\] 
is an $\mathbb{N}$-gradable $V$-module quotient of $W \cong \overline{M}_n(W(n)) /(\mathrm{ker} \, \Phi)$ by  $(\mathrm{ker} \, \Phi) /\overline{J}$, proving the first paragraph of the theorem.  

Now suppose that $V w \cap W(n)  = 0$ implies that $w = 0$ for $w \in W$.   Let $\bar{w} \in \bar{J}$ such that $\bar{w} = w + \mathcal{U}(\hat{V}) W_A$.  Then $\langle u', x w \rangle = 0$ for all $x \in \mathcal{U}(\hat{V})$ and $u' \in W(n)^*$.  Thus $V\cdot \Phi(\bar{w}) \cap W(n) = 0$, which implies that $\Phi (\bar{w}) = 0$.  Therefore $\overline{J} \subset \mathrm{ker} \, \Phi$, implying $\overline{J} = \mathrm{ker} \, \Phi$, proving the second paragraph of the theorem.
\epfv

\begin{rema}\label{J-remark} {\em Theorem \ref{last-thm} above gives some motivation and intuition about the subspace $J$ of $M_n(U)$ used to define $L_n(U)$.  In fact, from Theorem \ref{last-thm}, we see that $\overline{J} = J/\mathcal{U}(\hat{V})W_A$ is the maximal submodule of $\overline{M}_n(U)$ of the form $N/\mathcal{U}(\hat{V})W_A$ such that $N \cap U = 0$. } 
\end{rema}

We have the following corollary:

%
%

\begin{cor}\label{iso-cor}  Let $\mathcal{A}_{n,n-1}^{Res}$ denote the subcategory of objects $U$ in $\mathcal{A}_{n,n-1}$ that satisfy $\Omega_j(L_n(U)) = \bigoplus_{k=0}^j L_n(U)(k)$ for $j = n$ and $n-1$, and let $\mathcal{V}_n^{Res}$ denote the subcategory of objects $W$ in  $\mathcal{V}_n$ that satisfy: $\Omega_j(W) = \bigoplus_{k=0}^j W(k)$ for $j = n$ and $n-1$;
For any $w \in W$, $w = 0$ if and only if $Vw\cap W(n) = 0$;
and $W(n)$ has no nonzero $A_n(V)$-submodule that is an $A_{n-1}(V)$-module.  

Then the functors $\Omega_n/\Omega_{n-1}$ and $L_n$ are mutual inverses on the categories $\mathcal{V}_n^{Res}$ and $\mathcal{A}_{n,n-1}^{Res}$, respectively.  In particular, the categories $\mathcal{V}_n^{Res}$ and $\mathcal{A}^{Res}_{n,n-1}$ are equivalent. 

Furthermore, the subcategory of simple objects in $\mathcal{V}_n^{Res}$ is equivalent to the subcategory of simple objects in $\mathcal{A}_{n,n-1}^{Res}$. 
\end{cor}

\pf  By Theorem \ref{mainthm}, the functor $\Omega_n/\Omega_{n-1}$ takes objects in $\mathcal{V}_n^{Res}$ to objects in $\mathcal{A}_{n.n-1}^{Res}$, and $\Omega_n/\Omega_{n-1} \circ L_n$ is the identity on $\mathcal{A}_{n.n-1}^{Res}$.

By Theorem \ref{last-thm} and Remark \ref{J-remark}, the functor $L_n$ takes objects in $\mathcal{A}_{n.n-1}^{Res}$ to objects in $\mathcal{V}_n^{Res}$, and $L_n \circ \Omega_n/\Omega_{n-1}$ is the identity on $\mathcal{V}_n^{Res}$.

The correspondence on the subcategories of simple objects follows from Theorem \ref{simple-theorem}.  
\epfv

In particular, this illustrates that the categorical correspondence determined by the functors $\Omega_n/\Omega_{n-1}$ and $L_n$ restricts to subcategories of $\mathcal{V}_n$ and $\mathcal{A}_{n, n-1}$ that are too narrow to give a significant understanding of the relationship between indecomposable modules and higher level Zhu algebras as regards those modules that can be constructed via the functor $L_n$, which is the setting that motivated this paper in the first place.  That is, we are more interested in understanding the nature of the types of indecomposable $V$-modules that can be constructed from various classes of $A_n(V)$-modules through $L_n$, and in fact the functor $\Omega_n/\Omega_{n-1}$ is more useful in the indecomposable nonsimple setting as giving more information when it is {\it not} an inverse to $L_n$.  We study this issue further in \cite{BVY-Virasoro}, and in Section \ref{Virasoro-factor-through-example} below, we give an example to illustrate the types of indecomposable modules one can construct through the functor $L_n$ from indecomposable $A_n(V)$-modules. 

We also observe the following:

\begin{prop}\label{2n-prop} Let $U$ be an $A_n(V)$-module that does not factor through $A_{n-1}(V)$, and let $W = L_n(U)$.  Then
$\bigoplus_{k=0}^n W(k) \subset \Omega_n(W) \subset \bigoplus _{k=0}^{2n} W(k)$, and all singular vectors, $\Omega_0(W)$, must be contained in $\bigoplus _{k=0}^{n} W(k)$.
\end{prop}

\pf The first inclusion $\bigoplus_{k=0}^n W(k) \subset \Omega_n(W)$ is obvious.  Now suppose $w \in \Omega_n(W)$, and $w = w' + w''$ with $w' \in  \bigoplus _{k=0}^{2n} W(k)$ and $w'' \in  \bigoplus _{k=2n+1}^{\infty} W(k)$.  Then the $\mathbb{N}$-grading of $W$ implies that $w', w'' \in \Omega_n(W)$.  This and the fact that $U = W(n)$, implies that $\mathcal{U}(\hat{V}) w'' \cap U = 0$, since $ \Omega_n(W) \cap \bigoplus _{k=2n+1}^{\infty} W(k) = 0$.  Thus  $w'' \in J$, and $w = w'$ in $W = L_n(U) = M_n(U)/J$.  Therefore $\Omega_n(W) \subset \bigoplus _{k=0}^{2n} W(k)$. Furthermore, any singular vector $w$, i.e., any  $w \in \Omega_0(W)$,  must in fact then be contained in $\bigoplus _{k=0}^{n} W(k)$, otherwise if $w = w' + w''$ with $w' \in \bigoplus _{k=0}^{n} W(k)$ and $w'' \in \bigoplus_{k = n+1}^{\infty} W(k)$, then again $\mathcal{U}(\hat{V}) w'' \cap U = 0$, implying $w'' \in J$.  Therefore $\Omega_0(W) \subset \bigoplus _{k=0}^{n} W(k)$.
\epfv

\begin{rema} {\em Proposition \ref{2n-prop} and the fact that if $W = L_n(W(n))$ for $w \in W$, if $Vw = W$, then $V w \cap W(n) = 0$ if and only if $w = 0$ (or more generally, for any $w \in W$, then $\langle W(n)^*, w \rangle = 0$ if and only if $w = 0$), help to characterize the modules in $\mathcal{V}_n$ that are in the image of the functor $L_n$.   That is if $W = L_n(W(n))$, then $W$ must satisfy the following:

(i) $W$ is generated by $W(n)$;

(ii)  For $w \in W$, if $Vw = W$, then $V w \cap W(n) = 0$ if and only if $w = 0$ (or more generally, for any $w \in W$, then $\langle W(n)^*, w \rangle = 0$ if and only if $w = 0$);

(iii) $\Omega_n(W) \subset \bigoplus _{k=0}^{2n} W(k)$ and $\Omega_0(W) \subset \bigoplus _{k=0}^{n} W(k)$.}
\end{rema}

In particular, we see that the higher level Zhu algebras, i.e., $A_n(V)$ for $n \geq 1$, can be used to construct indecomposable nonsimple modules for $V$ with Jordan blocks for the $L(0)$ operator of sizes $k$ through $k+n$, if $A_n(V)$ does not decompose into a direct sum with $A_{n-1}(V)$; see Corollary \ref{mainthm-first-cor}.  We illustrate this below for $n = 1$ in Section \ref{Virasoro-factor-through-example}.

\section{Examples: Heisenberg and Virasoro vertex operator algebras}\label{examples-section}

In \cite{BVY-Heisenberg} and \cite{BVY-Virasoro}, we determine the level one Zhu algebra of $V$ in the two cases of when $V$ is the vertex operator algebra associated to the rank one Heisenberg algebra and when $V$ is the Virasoro vertex operator algebra, respectively.   In \cite{BVY-Heisenberg} and \cite{BVY-Virasoro}, we then give some results on the classification of modules for these vertex operator algebras using the structure of their level one Zhu algebras.  

Here we recall from \cite{FZ} and \cite{W} the level zero Zhu algebras and from \cite{BVY-Heisenberg} and \cite{BVY-Virasoro} the level one Zhu algebras for these vertex operator algebras.  We point out some distinctive features of these algebras which give interesting examples to illustrate aspects of Theorems \ref{mainthm} and \ref{simple-theorem}, as well as the other results of Section 3.  For instance, these examples illustrate how the nature of the surjection of the level one Zhu algebra onto the level zero Zhu algebra affects the structure of the $V$-modules that arise, in particular, how the nature of this surjection determines whether the higher level Zhu's algebras are necessary in order to detect indecomposable modules that have different size Jordan blocks with respect to $L(0)$ at higher degrees.

\subsection{The Heisenberg vertex operator algebra, and its level zero and level one Zhu algebras}\label{Heisenberg-section}

Following, for example \cite{LL}, we denote by $\mathfrak{h}$ a one-dimensional abelian Lie algebra spanned by $\alpha$ with a bilinear form $\langle \cdot, \cdot \rangle$ such that $\langle \alpha, \alpha \rangle = 1$, and by
\[
\hat{\mathfrak{h}} = \mathfrak{h}\otimes \C[t, t^{-1}] \oplus \C \mathbf{k}
\]
the affinization of $\mathfrak{h}$ with bracket relations
\[
[a(m), b(n)] = m\langle a, b\rangle\delta_{m+n, 0}\mathbf{k}, \;\;\; a, b \in \mathfrak{h},
\]
\[
[\mathbf{k}, a(m)] = 0,
\]
where we define $a(m) = a \otimes t^m$ for $m \in \mathbb{Z}$ and $a \in \mathfrak{h}$.

Set
\[
\hat{\mathfrak{h}}^{+} =\mathfrak{h} \otimes   t\C[t] \qquad \mbox{and} \qquad \hat{\mathfrak{h}}^{-} = \mathfrak{h} \otimes t^{-1}\C[t^{-1}].
\]
Then $\hat{\mathfrak{h}}^{+}$ and $\hat{\mathfrak{h}}^{-}$ are abelian subalgebras of $\hat{\mathfrak{h}}$.  Consider the induced $\hat{\mathfrak{h}}$-module given by 
\[
M(1) = \mathcal{U}(\hat{\mathfrak{h}})\otimes_{\mathcal{U}(\C[t]\otimes \mathfrak{h} \oplus \C \mathbf{k})} \C{\bf 1} \simeq S(\hat{\mathfrak{h}}^{-}) \qquad \mbox{(linearly)},
\]
where $\mathcal{U}(\cdot)$ and $S(\cdot)$ denote the universal enveloping algebra and symmetric algebra, respectively, $\mathfrak{h} \otimes \C[t]$ acts trivially on $\mathbb{C}\mathbf{1}$ and $\mathbf{k}$ acts as multiplication by $1$.   Then $M(1)$ is a vertex operator algebra, often called the {\it vertex operator algebra associated to the rank one Heisenberg algebra}, or simply the {\it rank one Heisenberg vertex operator algebra}, or the {\it one free boson vertex operator algebra} --- the Heisenberg Lie algebra in question being precisely $\hat{\mathfrak{h}} \diagdown \mathbb{C} \alpha(0)$.

There is in fact a one-parameter family of possible conformal elements for $M(1)$ that give the vertex operator algebra structure, namely $\omega_a = \frac{1}{2} \alpha(-1)^2{\bf 1} + a \alpha(-2){\bf 1}$ for $a \in \mathbb{C}$.  We distinguish these different vertex operator algebra structures on $M(1)$, by $M_a (1)$.

Any element of $M_a(1)$ can be expressed as a linear combination of elements of the form
\begin{equation}\label{generators-for-V}
\alpha(-k_1)\cdots \alpha(-k_j){\bf 1}, \quad \mbox{with} \quad  k_1 \geq \cdots \geq k_j \geq 1, \ \mbox{for $j \in \mathbb{N}$}.
\end{equation}

It is known that $M_a(1)$ is simple and has infinitely many nonisomorphic irreducible modules, which can be easily classified (see \cite{LL}). Furthermore, the indecomposable generalized modules have been completely determined, e.g., see \cite{M}.  In particular, we have

\begin{prop}[\cite{M}] Let $W$ be an indecomposable generalized $M_a(1)$-module.  Then as an $\hat{\mathfrak{h}}$-module  
\begin{equation}
W \cong M_a(1) \otimes \Omega(W) 
\end{equation}
where $\Omega(W) = \{w \in W \; | \; \alpha(n)w = 0 \mbox{ for all $n>0$} \}$ is the vacuum space.   
\end{prop}

\begin{rema} {\em  Note that in terms of the functors $\Omega_n$ for $n \in \mathbb{N}$, if $W$ is an indecomposable $M_a(1)$-module, then 
\[\mathbf{1} \otimes \Omega(W) = \Omega_0(W) .\]}
\end{rema}

We have the following level zero and level one Zhu algebras for $M_a(1)$ (cf. \cite{FZ}, \cite{BVY-Heisenberg}): \begin{equation}
A_0(M_a(1)) \cong \mathbb{C}[x,y]/(y-x^2) \cong \mathbb{C}[x]
\end{equation}
under the identification 
\begin{eqnarray}
\alpha(-1)\mathbf{1} + O_0(M_a(1)) &\longleftrightarrow& x + (p_0(x,y)), \label{x-for-A_0}\\
 \alpha(-1)^2\mathbf{1} + O_0(M_a(1)) &\longleftrightarrow& y + (p_0(x,y)), \label{y-for-A_0}
\end{eqnarray}  
where $p_0(x,y) = y - x^2$.

For the level one Zhu algebra, we have
\begin{eqnarray}
A_1(M_a(1)) &\cong& \mathbb{C}[x,y]/((y-x^2)(y-x^2-2)) \\
&\cong& \mathbb{C}[x,y]/(y-x^2)\oplus \mathbb{C}[x,y]/(y-x^2-2) \\
&\cong& \mathbb{C}[x] \oplus \mathbb{C}[x] \ \cong \ A_0(M_a(1)) \oplus \mathbb{C}[x]
\end{eqnarray}
under the identification 
\begin{eqnarray}
\alpha(-1)\mathbf{1} + O_1(M_a(1)) &\longleftrightarrow& x + (p_0(x,y)p_1(x,y)), \\ \alpha(-1)^2\mathbf{1} + O_1(M_a(1)) &\longleftrightarrow& y + (p_0(x,y)p_1(x,y)),
\end{eqnarray}  
where again $p_0(x,y) = y - x^2$, and in addition $p_1(x,y) = y - x^2 - 2$.

\begin{rema}\label{Heisenberg-remark} {\em In this case, since the ideals $I_0 = (p_0(x,y))$ and $I_1 = (p_1(x,y))$ are relatively prime, i.e. $I_0 + I_1 = \mathbb{C}[x,y]$, we have that the level one Zhu algebra is naturally isomorphic to a direct sum of $A_0(M_a(1)) \cong \mathbb{C}[x,y]/I_0$ and its direct sum complement which is isomorphic to $\mathbb{C}[x,y]/I_1$.  

Thus any indecomposable module $U$ for $A_1(M_a(1))$ will either be an indecomposable module for $A_0(M_a(1))$ or an indecomposable module for its direct sum complement, $\mathbb{C}[x,y]/I_1$.  That is we will either have that $U$ itself will factor through $A_0(M_a(1))$ or only the zero submodule of $U$ will factor through $A_0(M_a(1))$.

Therefore by Theorem \ref{mainthm}, any indecomposable $U$ module for $A_1(M_a(1))$ which does not factor through $A_0(M_a(1))$, will satisfy
\[U \cong \Omega_1/\Omega_0 (L_1(U)) ,\]
and in particular, any extra requirement that no nonzero submodule of $U$ factor through $A_0(M_a(1))$ is superfluous.  }
\end{rema}

We illustrate the points made in Remark \ref{Heisenberg-remark} explicitly below.

\subsubsection{Heisenberg Example:}\label{Heisenberg-example}

The indecomposable modules for $A_0(M_a(1)) \cong \mathbb{C}[x,y]/(y-x^2) \cong \mathbb{C}[x]$ are given by 
\begin{equation}\label{U-for-A_0}
U_0(\lambda, k) = \mathbb{C}[x,y]/((y-x^2), (x - \lambda)^k) \cong \mathbb{C}[x]/(x - \lambda)^k
\end{equation}
for $\lambda \in \C$ and $k \in \mathbb{Z}_+$, and 
\[
L_0(U_0(\lambda, k)) \cong M_a(1) \otimes_{\mathbb{C}} \Omega(\lambda, k),
\]
where $\Omega(\lambda,k)$ is a $k$-dimensional vacuum space such that $\alpha(0)$  
acts with Jordan form given by
\begin{equation}
\left[ \begin{array}{cccccc}
\lambda & 1 & 0 & \cdots  & 0 & 0\\
0 & \lambda & 1 & \cdots  & 0 & 0 \\
0 & 0 & \lambda & \cdots  & 0 & 0 \\
\vdots & \vdots & \vdots & \ddots & \vdots & \vdots\\
0 & 0 & 0 & \cdots  & \lambda & 1 \\
0 & 0 & 0 & \cdots  & 0 & \lambda
\end{array}
\right] .
\end{equation}
Note then that the zero mode of $\omega_a$ which is given by
\begin{equation}
L(0) = \sum_{m \in \mathbb{Z}_+} \alpha(-m) \alpha(m) + \frac{1}{2} \alpha(0)^2 - a \alpha(0) 
\end{equation}
acts on $\Omega(\lambda, k)$ such that the only eigenvalue is $\frac{1}{2} \lambda^2 - a \lambda$ (which is the lowest conformal weight of  $M_a(1) \otimes \Omega(\lambda, k)$) and $L(0) -( \frac{1}{2} \lambda^2 - a \lambda) Id_k$ with respect to a Jordan basis for $\alpha(0)$ acting on $\Omega(\lambda, k)$ is given by 
\begin{equation}
\left[ \begin{array}{ccccccc}
0 & \lambda-a & \frac{1}{2}- a  &0 &  \cdots  & 0 & 0 \\
0 & 0 & \lambda  - a  & \frac{1}{2} - a & \cdots  & 0 & 0  \\
0 & 0 & 0 & \lambda - a &  \cdots  & 0 & 0  \\
\vdots & \vdots & \vdots & \ddots & \vdots & \vdots \\
0 & 0 & 0 & 0& \cdots   & \lambda - a& \frac{1}{2} - a\\
0 & 0 & 0 & 0 &  \cdots  &0  & \lambda - a \\
0 & 0 & 0 & 0 & \cdots  & 0 & 0
\end{array}
\right] .
\end{equation}
Also note that $L(0)$ is diagonalizable if and only if:  (i) $k = 1$ which corresponds to the case when $M_a(1) \otimes \Omega(\lambda, k)$ is irreducible; (ii) $k = 2$ and $\lambda = a$; or (iii) $k>2$ and $\lambda = a = \frac{1}{2}$. 

These $M_a(1) \otimes \Omega(\lambda, k)$ exhaust all the indecomposable generalized $M_a(1)$-modules and the $\mathbb{N}$-grading of $M_a(1) \otimes \Omega(\lambda, k)$ is explicitly given by
\[ M_a(1) \otimes \Omega(\lambda, k) = \coprod_{m \in \mathbb{N}} M_a(1)_m \otimes \Omega( \lambda, k) \]
where $M_a(1)_m$ is the weight $m$ space of the vertex operator algebra $M_a(1)$ and thus $M_a(1)_m \otimes \Omega( \lambda, k)$ is the space of generalized eigenvectors of 
weight $m + \frac{1}{2}\lambda^2 - a \lambda$ with respect to $L(0)$.  Therefore, the generalized graded dimension of $M_a(1) \otimes \Omega(\lambda, k)$ is given by
\begin{eqnarray}
\mathrm{gdim}_q \, M_a(1) \otimes \Omega(\lambda, k) &=& q^{-1/24} \sum_{m \in \mathbb{N}} (k \,\mathrm{dim} \, M_a(1)_m )\, q^{m + \frac{1}{2} \lambda^2 - a \lambda}\\
&=&  q^{\frac{1}{2} \lambda^2 - a \lambda - 1/24} k \sum_{m \in \mathbb{N}} (\mathrm{dim} \, M_a(1)_m )\, q^{m} \nonumber  \\
&=&  q^{\frac{1}{2} \lambda^2 - a \lambda} k \, \eta(q)^{-1}  \nonumber
\end{eqnarray}
where $\eta(q)$ is the Dedekind $\eta$-function.

The indecomposable modules for $A_1(M_a(1)) \cong \mathbb{C}[x,y]/((y-x^2)(y - x^2 - 2)) \cong A_0(M_a(1)) \oplus \mathbb{C}[x]$ are given by the indecomposable modules $U_0(\lambda, k)$ for $A_0(M_a(1))$ as given in (\ref{U-for-A_0}) or by
\begin{equation}
U_1(\lambda, k) = \mathbb{C}[x,y]/((y-x^2 - 2), (x - \lambda)^k) \cong \mathbb{C}[x]/(x - \lambda)^k
\end{equation}
for $\lambda \in \C$ and $k \in \mathbb{Z}_+$, in which case, 
\[ U_1(\lambda, k) \cong \alpha(-1)\mathbf{1} \otimes \Omega(\lambda, k) \quad \mbox{and}  
\quad L_1(U_1(\lambda, k)) \cong M_a(1) \otimes \Omega(\lambda,k).\]
Therefore, there are no new modules obtained via inducing from a module for $A_1(M_a(1))$ versus from $A_0(M_a(1))$.  

If we allow for inducing by $L_1$ for any indecomposable $A_1(M_a(1))$  module (including those that factor through $A_0(M_a(1))$ so that $L_1(U)(0)$ might be zero), then the possible cases for $\Omega_1/\Omega_0(L_1(U))$ for $U$ an indecomposable $A_1(M_a(1))$-module are
\[ U = U_0(\lambda, k), \quad L_1(U) (0) = 0, \quad  \mbox{and} \quad \Omega_1/\Omega_0(L_1(U)) \cong U_1(\lambda, k) \ncong U\]
or 
\[ U = U_1(\lambda, k)  \quad L_1(U) (0) = \Omega(\lambda, k) \neq  0, \quad \mbox{and} \quad \Omega_1/\Omega_0(L_1(U)) \cong U.\]
Note however that in the case of $U = U_0(\lambda,k)$, the $M_a(1)$-module $L_1(U)$ is in fact $M_a(1) \otimes \Omega(\lambda, k)$ but the grading as an $\mathbb{N}$-gradable module is shifted up one.   Thus by regrading to obtain an $\mathbb{N}$-gradable $M_a(1))$ module in the sense of Definition \ref{N-gradable-definition}, this module is again just $M_a(1) \otimes \Omega(\lambda, k)$, and the level one Zhu algebra gives no new information about the indecomposable $M_a(1)$-modules not already given by the level zero Zhu algebra.

\subsection{Virasoro vertex operator algebras, and their level zero and level one Zhu algebras}

Let $\mathcal{L}$ be the Virasoro algebra with central charge $\mathbf{c}$, that is, $\mathcal{L}$ is the vector space with basis $\{\bar{L}_n \,|\, n\in\Z\}\cup \{\mathbf{c}\}$ with bracket relations 
\begin{align*}
[\bar{L}_m,\bar{L}_n]=(m-n)\bar{L}_{m+n}+\frac{m^3-m}{12} \delta_{m+n,0} \, \textbf{c},\quad\quad  [\textbf{c},\bar{L}_m]=0
\end{align*}  
for $m,n\in \Z$.  Here we use a bar over the Virasoro generators to distinguish between these Virasoro elements and the functor $L_n$ defined earlier.  

Let $\mathcal{L}^{\geq 0}$ be the Lie subalgebra with basis $\{ \bar{L}_n  \,|\, n\geq 0 \} \cup \{\mathbf{c}\}$, and let $\mathbb{C}_{c,h}$ be the $1$-dimensional $\mathcal{L}^{\geq 0}$-module where $\mathbf{c}$ acts as $c$ for some $c\in \mathbb{C}$, $\bar{L}_0$ acts as $h$ for some $h\in \mathbb{C}$, and $\bar{L}_n$ acts trivially for $n\geq 1$. Form the induced $\mathcal{L}$-module
\[
M(c,h)= \mathcal{U}(\mathcal{L})\otimes_{\mathcal{L}^{\geq 0}} \mathbb{C}_{c,h} .
\]
We shall write $L(n)$ for the operator on a Virasoro module corresponding to $\bar{L}_n$, and $\mathbf{1}_{c,h} = 1 \in \mathbb{C}_{c,h}$.  Then 
\[ V_{Vir}(c,0)= M(c,0)/\langle L(-1)\mathbf{1}_{c,0}\rangle
\] 
has a natural vertex operator algebra structure with vacuum vector $1=\mathbf{1}_{c,0}$, and conformal element $\omega=L(-2)\mathbf{1}_{c,0}$, satisfying  $Y(\omega,x)= \sum_{n\in\Z } L(n)x^{-n-2}$.  In addition, for each $h\in\C$, we have that $M(c,h)$ is an ordinary $V_{Vir}(c,0)$-module with $\mathbb{N}$-gradation 
\[M(c,h)=\coprod_{k\in \N} M(c,h)_k\] 
where $M(c,h)_k$ is the $L(0)$-eigenspace with eigenvalue $h + k$.  We say that  $M(c,h)_k$ has degree $k$ and weight $h+k$.   

We now fix $c \in \mathbb{C}$, and denote by $V_{Vir}$, the vertex operator algebra $V_{Vir}(c,0)$.  

It was shown in \cite{W} that
\begin{equation}
A_0(V_{Vir}) \cong \mathbb{C}[x,y]/(y-x^2 - 2x) \cong \mathbb{C}[x]
\end{equation}
under the identification 
\begin{eqnarray}
L(-2)\mathbf{1} + O_0(V_{Vir}) &\longleftrightarrow& x + (q_0(x,y)), \\ 
L(-2)^2\mathbf{1} + O_0(V_{Vir}) &\longleftrightarrow& y + (q_0(x,y)),
\end{eqnarray}  
where $q_0(x,y) = y - x^2 - 2x$.

In addition, there is a bijection between isomorphism classes of irreducible $\mathbb{N}$-gradable $V_{Vir}$-modules and irreducible $\mathbb{C}[x]$-modules given by $L(c,\lambda)
\longleftrightarrow \mathbb{C}[x]/(x - \lambda)$ where $T(c,\lambda)$ is the largest proper submodule of $M(c,\lambda)$, and $L(c,\lambda) = M(c,\lambda)/T(c,\lambda)$.  It was proved in \cite{W} that if $c$ is not of the form 
\begin{equation}\label{cpq}
c_{p,q} = 1 - 6\frac{(p-q)^2}{pq} \quad \mbox{for any $p,q \in \{2,3,4,\dots\}$ with $p$ and $q$ relatively prime,}
\end{equation}
then $T(c, \lambda) = \langle L(-1) \mathbf{1}_{c,0} \rangle$, implying $V_{Vir}(c,0) = L(c,0)$, is simple as a vertex operator algebra.  Furthermore, it then follows from the fact $A_0(V_{Vir}) \cong \mathbb{C}[x]$ that in this case of $c \neq c_{p,q}$, $V_{Vir}(c,0)$ is irrational.   It was also shown in \cite{W}  that for $c = c_{p,q}$ that $T(c, \lambda) \neq \langle L(-1) \mathbf{1}_{c,0}$ and thus in this case $V_{Vir}(c_{p,q}, 0)$ is not simple as a vertex operator algebra.  

As for the level one Zhu algebra, in \cite{BVY-Virasoro}, we show that 
\begin{eqnarray}
A_1(V_{Vir}) &\cong& \mathbb{C}[x,y]/((y-x^2-2x)(y-x^2-6x + 4)) \\
&\cong& \mathbb{C}[\tilde{x},\tilde{y}]/(\tilde{x} \tilde{y})
\end{eqnarray}
under the identification 
\begin{eqnarray}
L(-2)\mathbf{1} + O_1(V_{Vir}) &\longleftrightarrow& x + (q_0(x,y)q_1(x,y)),  \label{x}\\ 
L(-2)^2\mathbf{1} + O_1(V_{Vir}) &\longleftrightarrow& y + (q_0(x,y)q_1(x,y)), \label{y}
\end{eqnarray}  
where $q_0(x,y) = y - x^2 -2x = \tilde{x}$ and in addition $q_1(x,y) = y - x^2 -6x + 4 = \tilde{y}$.

\begin{rema} {\em In this case, the ideals $J_0 = (q_0(x,y))$ and $J_1 = (q_1(x,y))$ are not relatively prime since $q_0(x,y) - q_1(x,y) = 4x - 4$, i.e., $q_0$ and $q_1$ have nontrivial intersection at $x=1$,  and thus $J_0 + J_1 \neq \mathbb{C}[x,y]$.  Therefore the level zero Zhu algebra is not isomorphic to a direct summand of the level one Zhu algebra.  This results, as we shall see below, in several interesting examples illustrating the subtleties  of the relationship between modules for  the higher level Zhu algebras for $V = V_{Vir}$ and $\mathbb{N}$-gradable $V_{Vir}$-modules. }  
\end{rema}

\begin{rema} {\em The result we obtain in  \cite{BVY-Virasoro} for $A_1(V_{Vir})$ differs from that presented in \cite{V} as follows: In the notation of \cite{V}, letting $L = L(0)$ and $A = L(-1)L(1)$, then acting on a primary vector (i.e. a vector $w$ for which $L(n) w = 0$ if $n >0$), we have that $(y - x^2 -2x)(y-x^2 -6x + 4)$ acts (using the zero mode action as computed below in (\ref{L(0)-zero-mode}) and (\ref{L(-2)-squared-acting})) as $4(A^2 - 2AL + 2A)$.  This implies that for such vectors $A_1(V_{Vir})$  acts as $\mathbb{C}[A, L]/(A^2 - 2LA +2A)$ and this algebra is almost the level one Zhu algebra for $V_{Vir}$ given in \cite{V} but still differs by a minus sign.  Thus with the minus sign typo corrected, $A_1(V_{Vir})$ acts equivalently to the algebra given in \cite{V} on $\Omega_1(W)$ for an $\mathbb{N}$-gradable $V$-module $W$, but will not in general act the same on $\Omega_n(W)$ for $n>1$.  Since one of the important aspects of the information that the higher level Zhu algebras give is contained in the action of $A_n(V)$ versus the action of $A_{n-1}(V)$ through the natural epimorphism from $A_n(V)$ to $A_{n-1}(V)$, it is essential to have each $A_n(V)$ realized as its full algebra $V/O_n(V)$ rather than as a realization of zero modes acting on $\Omega_n(V)$ so as to be able to compare, e.g. the action of $A_2(V_{Vir})$ on an $A_2(V_{Vir})$-module such as $\Omega_2(W)$ versus the action of $A_1(V_{Vir})$ on the same module. }
\end{rema}

\subsubsection{Virasoro Example 1:  $\Omega_n/\Omega_{n-1}(L_n(U)) \ncong U$ since $U$ has a nonzero proper submodule that factors through $A_{n-1}(V)$}\label{Virasoro-factor-through-example}

We now give a family of examples to illustrate that there are nontrivial instances of an indecomposable $A_1(V_{Vir})$-module $U$ that does not factor through $A_0(V_{Vir})$, but a nontrivial submodule does, and that in this case we have
\[\Omega_1/\Omega_0 (L_1(U)) \ncong U.\] 
In addition, we give other aspects of the structure of this family of indecomposable $V_{Vir}$-modules, $\{W_k\}_{k \in \mathbb{Z}_+}$.  For instance, each $W_k$ for $k \in \mathbb{Z}_+$, has a Jordan block decomposition for $L(0)$ with Jordan blocks of size $k$ in the weight zero, and weight greater than one subspaces, and size $k+1$ in the weight one and greater subspaces. 

For $k \in \mathbb{Z}_+$, consider the non simple, indecomposable $A_1(V_{Vir})$-module given by 
\[
U = \C[x, y]/((y-x^2-2x)^{k+1}, (y-x^2-6x+4)).
\]
Clearly $U$ is not a module for $A_0(V_{Vir}) \cong \mathbb{C}[x,y]/(y - x^2 - 2x)$, and
\[
U \cong \C[x]/(x-1)^{k+1}.
\]
Let $w$ be a singular vector for the Virasoro algebra (i.e., $L(n) w = 0$ if $n>0$) such that
\[
L(0)^kw \neq 0; \;\;\; L(0)^{k+1}w = 0.
\]
Set
\[
U' = \mathrm{span}_\mathbb{C} \{L(-1)L(0)^iw\;|\; i = 0, \dots, k\}.
\]

We have the following lemma:
\begin{lemma}
As $A_1(V_{Vir})$-modules, $U \cong U'$
under the homomorphism 
\begin{eqnarray}
f : \ \ \ \  \ \ \ U &\longrightarrow&  U'\\
\overline{(x-1)^{i}} &\mapsto& L(-1)L(0)^iw, \nonumber
\end{eqnarray}
for $i = 0, \dots, k$, and
where $\overline{(x-1)^{i}}$ is the image of $(x-1)^i$ in $U$ under the canonical projection.
\end{lemma}
\pf Clearly the map $f$ is surjective, and since $U$ and $U'$ have the same dimension, $f$ is also injective. We need to show $f$ is an algebra homomorphism.

Denote by $[x]$ and $[y]$ be the image of $x$ and $y$ in $A_1(V_{Vir})$ under the canonical projection and identification given by (\ref{x}) and (\ref{y}).  Then $[x] = [L(-2){\bf 1}]$ and $[y] = [L(-2)^2{\bf 1}]$. Thus $[x]$ acts on the modules via
\begin{equation}\label{L(0)-zero-mode}
o(L(-2){\bf 1}) = L(0). 
\end{equation}
To determine the action of $[y]$, recall the normal ordering notation 
\[{}^\circ_\circ  Y(u,x)Y(v,x)   {}^\circ_\circ = \left(\sum_{m<0} u_n x^{-n-1}\right) Y(v,x) + Y(v,x) \sum_{m \geq 0} u_n x^{-n-1}.\] 
Then 
\begin{eqnarray}\label{L(-2)-squared-acting}
o(L(-2)^2{\bf 1}) &=& (L(-2)^2{\bf 1})_{\mathrm{wt} \, L(-2)^2{\bf 1} - 1} \ = \ (L(-2)^2{\bf 1})_3 \nn
&=& \mathrm{Res}_x x^3 \,  {}^\circ_\circ Y(L(-2) \mathbf{1},x) Y(L(-2) \mathbf{1},x)  {}^\circ_\circ \nn
&=& \sum_{m<0} L(m-1)L(-m+1) + \sum_{m\geq 0}L(-m+1) L(m-1) \nn
&=& L(1)L(-1) +  L(0)^2 + L(-1)L(1) + 2 \sum_{i \geq 2}L(-i)L(i) \nn
&=& [L(1), L(-1)] + L(0)^2 + 2 \sum_{i \geq 1}L(-i)L(i) \nn
&=& 2L(0) + L(0)^2 + 2 \sum_{i \geq 1}L(-i)L(i) .
\end{eqnarray}
Thus we have
\begin{eqnarray*}
f([x]\cdot \overline{(x-1)^i}) &=& f(\overline{(x-1)^{i+1}} + \overline{(x-1)^{i}}) \\
&=& f(\overline{(x-1)^{i+1}}) + f(\overline{(x-1)^{i}})\\
&=& L(-1)L(0)^{i+1}w + L(-1)L(0)^iw\\
&=& L(0)(L(-1)L(0)^iw)\\
&=& o(L(-2){\bf 1})(L(-1)L(0)^iw)\\
&=& [x]\cdot f( \overline{(x-1)^i}),
\end{eqnarray*}
and 
\begin{eqnarray*}
f([y]\cdot \overline{(x-1)^i}) &=& f([x^2+6x-4]\cdot \overline{(x-1)^i})\\
&=& f([(x-1)^2 + 8(x-1) + 3] \cdot \overline{(x-1)^i})\\
&=& f(\overline{(x-1)^{i+2}}) + 8f(\overline{(x-1)^{i+1}} ) + 3f(\overline{(x-1)^i})\\
&=& L(-1)L(0)^{i+2}w + 8L(-1)L(0)^{i+1}w + 3 L(-1)L(0)^iw,
\end{eqnarray*}
while from (\ref{L(-2)-squared-acting}), we have
\begin{eqnarray*}
\lefteqn{[y]\cdot f( \overline{(x-1)^i}) =}\\
&=& \left( 2L(0) + L(0)^2 + 2 \sum_{j \geq 1}L(-j)L(j) \right)\cdot L(-1)L(0)^iw\\
&=& L(-1)L(0)^{i+2}w + 8L(-1)L(0)^{i+1}w +  3L(-1)L(0)^iw\\
&=& f([y]\cdot \overline{(x-1)^i}).
\end{eqnarray*}
Therefore $f$ is a $A_1(V_{Vir})$-module homomorphism, proving the lemma. \epfv

Now we see that in this case 
\[
\Omega_1(L_1(U)) = L_1(U)(0) \oplus L_1(U)(1)  \oplus \mbox{higher degree terms}
\]
and
\[
\Omega_0(L_1(U)) = L_1(U)(0) \oplus  \left( \mathrm{span}_\mathbb{C} \{L(-1)L(0)^k w\} + J \right) .
\]
Since $L(-1)L(0)^kw \in U$, by definition of $J$, we have that $L(-1)L(0)^kw \notin J$. In particular,  
\begin{eqnarray*}
\Omega_1/\Omega_0(L_1(U)) &\cong& U / \mathrm{span}_\mathbb{C} \{L(-1)L(0)^kw\} \oplus  \mbox{higher degree terms}\\
&\ncong& U,
\end{eqnarray*}
giving a counter example to Theorem 4.2 in \cite{DLM} and showing the necessity of the condition in our Theorem \ref{mainthm} that no nontrivial submodule of $U$ factors through $A_{n-1}(V)$.  That is this case illustrates that since the submodule $\mathrm{span}_\mathbb{C}\{L(-1)L(0)^kw\} \subsetneq U$ is an $A_0(V_{Vir})$-module, the added condition in our Theorem \ref{mainthm} in comparison to Theorem 4.2 of \cite{DLM} is indeed necessary for the statement to hold.  

Note that in general here, $J$ and the nature of the higher degree terms in $\Omega_1(L_1(U))$, will depend on the central charge $c$ of $V_{Vir} = V_{Vir}(c,0)$.  In general,
\begin{eqnarray*}
(L_1(U)) (j) &=& \mathrm{span}_\mathbb{C} \{L(-s_1) \cdots L(-s_r)L(-1) L(0)^i w \; | \; i = 0, \dots, k, \  r \in \mathbb{N}, \\
& & \qquad s_1 \geq s_2 \geq\dots \geq s_r \geq 1, \ s_1 + \cdots + s_r + 1 = j \} \\
& & \quad \oplus \mathrm{span}_\mathbb{C} \{L(-s_1) \cdots L(-s_r) L(0)^i w \; | \; i = 1, \dots, k, \  r \in \mathbb{N}, \\
& & \qquad s_1 \geq s_2 \geq\dots \geq s_r > 1, \ s_1 + \cdots + s_r  = j \}  \ \mathrm{mod} \ J
\end{eqnarray*}
for $j \in \mathbb{N}$, where the first direct summand, which occurs if $j \neq 0$, consists of Jordan blocks for $L(0)$ of size $k+1$ modulo the subspace $J$, and the second direct summand, which occurs if $j \neq 1$, consists of Jordan blocks of size $k$ modulo the subspace $J$.  That is, the degree zero space will have a Jordan block of size $k$, the degree one space will have Jordan block of size $k+1$, and the higher degree subspaces will potentially have Jordan blocks of size $k+1$ and lower with the size depending on the central charge $c \in \C$.



For instance, if $k = 1$, then $U = \mathrm{span}_\C \{ L(-1)w, L(-1)L(0)w \}$, and we have
\begin{eqnarray}
L_1(U)(0) &=& \mathrm{span}_\C \{ L(0)w \},  \qquad  \qquad  \qquad \quad \ \  \mbox{ 1 Jordan block of size 1}\\
L_1(U)(1) &=& \mathrm{span}_\C \{ L(-1)w, L(-1)L(0)w \}, \quad  \mbox{ 1 Jordan block of size 2} .
\end{eqnarray}

\subsubsection{Virasoro Example 2:  $L_n(\Omega_n/\Omega_{n-1}(W)) \ncong W$ since $W \neq \langle W(n) \rangle$ even if $W$ is simple}\label{Virasoro-n-generated-example}

Now let $c \neq 0$, and let $L(c,0)$ be the unique simple minimal vertex operator algebra with central charge $c$, up to isomorphism, i.e. $L(c,0)$ is isomorphic to the quotient of $V_{Vir} = V_{Vir}(c,0)$ by its largest proper ideal $T(c,0)$.  

Let $V = W = L(c,0)$.  Then $W(0) = \mathbb{C} \mathbf{1}$, with $\mathbf{1} = \mathbf{1}_{(c,0)}$, and $W(1) = 0$.  Thus in this case 
\[\Omega_0(W) = \mathbb{C} \mathbf{1}= W(0) = W(0) \oplus W(1) = \Omega_1(W).\]  Therefore 
\[L_1(\Omega_1/\Omega_0 (W)) = 0 \ncong W.\]

\subsubsection{Virasoro Example 3: $L_n(\Omega_n/\Omega_{n-1}(W)) \ncong W$ since $\Omega_n/\Omega_{n-1}(W) \neq  W(n)$, for $W$ indecomposable but nonsimple, and $V$ simple}

Now let $V = L(c,0)$  and $c \neq c_{p,q}$ for $c_{p,q}$ defined in  (\ref{cpq}).  In this case $V= L(c,0) = V_{Vir} (c,0)$ and $V$ is a simple vertex operator algebra.  Let  $W = M(c,0)$ which is not a simple $V$-module but is indecomposable.  In this case, since the quotient of $M(c,0)$ by $\langle L(-1) \rangle$ is simple, we have
\[\Omega_0(W) = \mathrm{span}_\mathbb{C} \{ \mathbf{1}, L(-1)\mathbf{1}\} = W(0) \oplus W(1) \]
\[\Omega_1(W) = \mathrm{span}_\mathbb{C} \{ \mathbf{1}, L(-1)\mathbf{1}, L(-1)^2\mathbf{1}\}  = W(0) \oplus W(1) \oplus \mathbb{C} L(-1)^2 \mathbf{1}.\]
Thus
\[\Omega_1/\Omega_0(W) \cong \mathbb{C}  L(-1)^2\mathbf{1}, \]
which does not factor through $A_0(V)$ since $y - x^2 - 2x$ acts as $2L(-1)L(1)$ on  $\mathbb{C}  L(-1)^2\mathbf{1}$, which is nontrivial.  Then $L_1(\mathbb{C}  L(-1)^2\mathbf{1}) \cong \langle L(-1)\mathbf{1} \rangle$ since $L(-1) \mathbf{1} = \frac{1}{2} L(1) L(-1)^2 \mathbf{1}$ and this spans the zero degree space of $L_1(\mathbb{C}  L(-1)^2\mathbf{1})$.

Therefore we have 
\[L_1(\Omega_1/\Omega_0(W)) = L_1(\mathbb{C}  L(-1)^2\mathbf{1}) \cong \langle L(-1)\mathbf{1} \rangle  \ncong M(c,0) = W.\]

Also note that $\Pi_1(W) = \mathbb{C} L(-1) \mathbf{1}$ which as an $A_1(V)$-module does factor through $A_0(V)$, and $L_0( \Pi_1(W)) \cong  \langle L(-1)\mathbf{1} \rangle$.

\end{document}